\documentclass[12pt]{article}

\usepackage{amssymb,amsthm,amsmath}

\title{On Weak Tail Domination of Random Vectors}
\author{Rafa{\l} Lata{\l}a
\thanks{Research partially supported by MEiN Grant 1 PO3A 012 29.}}
\date{}

\newtheorem{thm}{Theorem}

\newtheorem{prop}{Proposition}
\newtheorem{dfn}{Definition}
\newtheorem{rem}{Remark}

\def\Ex{{\mathbf E}}
\def\Pr{{\mathbf P}}
\def\ve{\varepsilon}
\def\conv{\mathrm{conv}}
\def\cl{\mathrm{cl}}

\begin{document}

\maketitle

\begin{abstract}
Motivated by a question of Krzysztof Oleszkiewicz we study a
notion of weak tail domination of random vectors.
We show that if the dominating random variable is sufficiently
regular weak tail domination implies strong tail domination.
In particular positive answer to Oleszkiewicz question would
follow from the so-called Bernoulli conjecture.
\end{abstract}

{\bf Introduction.}
This note is motivated by the
following problem about Rademacher series, posed by Krzysztof Oleszkiewicz
(private comunication):

{\bf Problem.} {\em Suppose that 
$(\ve_i)$ is a Rademacher sequence (i.e. sequence of independent symmetric 
$\pm 1$ r.v.'s) and $x_i,y_i$ are vectors in some Banach space $F$
such that the series $\sum_{i}x_{i}\ve_{i}$ and $\sum_{i}y_{i}\ve_{i}$
are a.s. convergent and
\[
\forall_{x^*\in F^*}\forall_{t>0}\
\Pr\Big(\Big|x^*\Big(\sum_{i}x_{i}\ve_{i}\Big)\Big|\geq t\Big)\leq
\Pr\Big(\Big|x^*\Big(\sum_i y_i\ve_i\Big)\Big|\geq t\Big).
\]
Does it imply that
\[
  \Ex\Big\|\sum_i x_i\ve_i\Big\|\leq L\Ex\Big\|\sum_i y_i\ve_i\Big\|,
\]
for some universal constant $L<\infty$?}

Motivated by the above question we introduce a notion
of weak tail domination of random vectors. We show that
if the dominating vector has a regular distribution
(including Gaussian case), 
weak tail domination yields
strong tail domination (Theorem \ref{str_tails}). In particular 
Oleszkiewicz question has positive answer provided
that the so-called Bernoulli Conjecture holds true. 
Finally we show that in general weak tail domination does not
yield comparison of means or medians of norms even if the
distribution of dominated vector is Gaussian.

\begin{dfn}
Let $X$ and $Y$ be random vectors with values in some Banach space
$F$. We say that {\em tails of $Y$ are weakly dominated
by tails of $X$} and denote it by $Y\prec_{\omega} X$ if
\[
  \Pr(|x^{*}(Y)|\geq t)\leq \Pr(|x^*(X)|\geq t)
 \mbox{ for all } x^*\in F^*, t>0.
\]
\end{dfn}

The following regularity property of random vectors will give us 
a tool to pass from weak to strong comparison.

\begin{dfn}
\label{K_reg}
We say that a random vector $X$ with values in  $F$ is
$K$-regular for some $K<\infty$ if there exists a sequence
$(x_n^*)\subset F^*$ such that
\[
\|x_n^*(X)\|_{\log(n+2)}=
(\Ex|x_n^*(X)|^{\log (n+2)})^{1/\log(n+2)}\leq K\Ex\|X\|
\mbox{ for } n=1,2,\ldots.
\]
and
\[
B_{F^*}=\{x^*\in F^*\colon \|x^*\|\leq 1\}
\subset
\cl_{X}(\conv\{\pm x_n^*\colon n\geq 1\}),
\]
where for $A\subset F^*$,
$\cl_X(A)$ denotes the closure of $A$ with respect to the $L^2$ distance
$d_X(x^*,y^*):=(\Ex|x^*(X)-y^*(X)|^2)^{1/2}$.
\end{dfn}

\begin{prop}
\label{str_exp}
If $X$ is $K$-regular and $Y\prec_{\omega} X$, then $\Ex\|Y\|\leq 20K\Ex\|X\|$.
\end{prop}

\begin{proof}
Let $x_n^*$ be as in  Definition \ref{K_reg}. We have
for any $t>0$,
\begin{align*}
\Pr\Big(\sup_{n\geq 1}|x_n^*(Y)|\geq t\Big)
&\leq \sum_{n\geq 1}\Pr(|x_n^*(Y)|\geq t)
\leq \sum_{n\geq 1} t^{-\log(n+2)} \Ex|x_n^*(Y)|^{\log(n+2)}
\\
\leq&
 \sum_{n\geq 1} t^{-\log(n+2)} \Ex|x_n^*(X)|^{\log(n+2)}
\leq \sum_{n\geq 1} \Big(\frac{K\Ex\|X\|}{t}\Big)^{\log(n+2)}.
\end{align*}
Notice that 
$d_Y(x^*,y^*)\leq d_X(x^*,y^*)$, hence $B_F^*$ is contained 
also in the closure of absolute convex of $\pm x_n^*$ with 
respect to $d_Y$ metric and thus
\begin{align*}
\Ex \|Y\|
&\leq \Ex\sup_{n\geq 1}|x_n^*(Y)|\leq
K\Ex\|X\|\Big(e^2+\int_{e^2}^{\infty}
\Pr\Big(\sup_{n\geq 1}|x_n^*(Y)|\geq tK\Ex\|X\|\Big)dt\Big)
\\
&\leq K\Ex\|X\|\Big(e^2+\sum_{n=1}^{\infty}\int_{e^2}^{\infty}t^{-\log(n+2)}dt
\Big)
\leq 20K\Ex\|X\|.
\end{align*}
\end{proof}

\begin{thm}
\label{str_tails}
Let $X_1,X_2,\ldots$ be independent copies of symmetric random vector $X$. 
Suppose that there exist constants $K<\infty$ and $\alpha,\beta>0$ such that
for all $n=1,2,\ldots$\\
i) random vector $(X_1,\ldots,X_n)$ with values in $l^n_{\infty}(F)$ is
$K$-regular,\\
ii) $\Pr(\max_{i\leq n}\|X_i\|\geq \alpha \Ex\max_{i\leq n}\|X_i\|)\geq \beta$.\\
Then for any  random vector $Y$ such that $Y\prec_{\omega} X$ we have
\[
\Pr(\|Y\|\geq t)\leq \frac{2}{\beta}\Pr\Big(\|X\|\geq \frac{\alpha t}{80 K}\Big).
\]
\end{thm}

The main idea how to derive comparison of tails from comparison of means 
is not new - it goes back at least to
the paper of de la Pe\~na, Montgomery-Smith and Szulga \cite{dPMSS}.

\begin{proof}
We may obvoiusly assume that $Y$ is symmetric, by $Y_1,Y_2,\ldots$ we will denote
independent copies of $Y$.
Let $n\geq 2$ be such that 
\[
\frac{2}{n}\geq \Pr(\|Y\|\geq t)\geq \frac{1}{n}.
\]
Then $\Pr(\max_{i\leq n}\|Y_i\|\geq t)\geq 1-(1-1/n)^n\geq 1/2$, hence
$\Ex\max_{i\leq n}\|Y_i\|\geq t/2$. Let $\eta$ be r.v. independent
of $(Y_i)$ such that
$\Pr(\eta=1)=\Pr(\eta=0)=1/2$, then by Theorem 3.2.1 of \cite{KW},
$\eta(Y_1,\ldots,Y_n)\prec_{\omega} (X_1,\ldots,X_n)$, where both variables
are considered as random vectors in $l^n_{\infty}(F)$. By Proposition
\ref{str_exp},
\begin{align*}
\frac{t}{4}&\leq \Ex\max_{i\leq n}\|\eta Y_i\|
=\Ex \|\eta(Y_1,\ldots,Y_n)\|_{l_{\infty}^n(F)}\leq 
20K\Ex \|(X_1,\ldots,X_n)\|_{l_{\infty}^n(F)}
\\
&=
20K\Ex\max_{i\leq n}\|X_i\|.
\end{align*}
Property ii) yields 
\[
\beta\leq \Pr\Big(\max_{i\leq n}\|X_i\|\geq \frac{\alpha t}{80K}\Big)
\leq n\Pr\Big(\|X\|\geq \frac{\alpha t}{80K}\Big),
\]
so $\Pr(\|X\|\geq \alpha t/(80K))\geq \beta/n\geq \beta \Pr(\|Y\|\geq t)/2.$
\end{proof}

\begin{rem}
\label{prop_2}
By the Paley-Zygmund inequality (cf.\ \cite{KW}, Lemma 0.2.1), 
the comparison of first and second moments of maxima,
\begin{equation}
\label{comp_mom}
\Ex\max_{i\leq n}\|X_i\|^2\leq C(\Ex\max_{i\leq n}\|X_i\|)^2
\end{equation}
implies 
property ii) of previous
theorem with $\alpha=1/2$ and $\beta=1/(4C)$.
\end{rem}

\begin{rem}
Both Proposition \ref{str_exp} and Theorem \ref{str_tails}
hold (with constants depending on $C_1$ and $C_2$) 
if we replace the condition $Y\prec_{\omega} X$ by the condition
\begin{equation}
\label{weak_tails2}
\Pr(|x^*(Y)|\geq t)\leq C_1\Pr(|x^*(X)|\geq t/C_2)
\mbox{ for all } x^*\in F^*, t>0.
\end{equation}
\end{rem}

Indeed, if $\eta$ is a random variable independent of $Y$ with
$\Pr(\eta=1)=1-\Pr(\eta=0)=1/C_1$, then condition (\ref{weak_tails2})
implies $\eta Y/C_2 \prec_{\omega} X$.

\medskip

Let us give few examples of random vectors satisfying the assumptions
of Theorem 1.

\medskip

1. Any centered Gaussian vector on a separable Banach space is
$L$-regular with universal $L$. This is a consequence of majorizing
measure theorem (cf.\cite{Ta2} and  \cite{Ta1}, Theorem 2.1.8).
Since a product of Gaussian measures is again Gaussian, property i)
holds with $K=L$. Moments of Gaussian vectors are comparable
so by Remark \ref{prop_2} also property ii) holds with $\alpha=1/2$
and universal $\beta$.

\medskip

2. Let $(\eta_i)$  be a sequence of independent symmetric real r.v.'s
 with logarithmically concave tails satisfying $\Delta_2$ condition and
$v_i\in F$ be such that $X=\sum_{i}v_i\eta_i$ is a.s. convergent.
Then  $X$ is $K$-regular with constant $K$ depending only on
$\Delta_2$ constant (\cite{La}, Theorem 3). Random variable $(X_1,\ldots,X_n)$
has an analogous series representation in $l^n_\infty(F)$, so property
i) holds. It can be also checked that (\ref{comp_mom}) is satisfied
with universal $C$.

\medskip

3. Positive answer to Bernoulli Conjecture (\cite{Ta1}, Chapter 4)
would imply the $L$-regularity of Rademacher series. Since (\ref{comp_mom})
holds for $X$ being a Rademacher sum with vector coefficients, Theorem
1 would give positive answer to Oleszkiewicz question.
 
\medskip

We conclude with an example showing that weak tail domination does not yield 
any comparison
of strong parameters even if the dominated vector has Gaussian
distribution. 

\medskip

{\bf Example.}
Let $F=l^n_{\infty}$, $Y=\sum_{i=1}^n g_ie_i$ and
$X=9(|g_1|+1)\sum_{i=1}^n\eta_ie_i$,
where $g_i$ are i.i.d. ${\cal N}(0,1)$ and
$\eta_i$ are i.i.d.\ r.v.'s with uniform distribution on $[-1,1]$, independent
of $g_1$. 

To show that  tails of $Y$ are weakly dominated by tails of $X$
it is enough to check that
\begin{equation}
\label{example_1}
\Pr(|\langle u,Y\rangle|\geq t)\leq \Pr(|\langle u,X\rangle|\geq t)
\mbox{ for } u\in S^{n-1},\ t\geq 0.
\end{equation}
Let us fix $u\in S^{n-1}$. For any $t>0$ we have
\[
\Pr(|\langle u,Y\rangle|\geq t)=\Pr(|g_1|\geq t).
\]
By the Paley-Zygmund inequality,
\begin{align*}
\Pr\Big(\Big|\sum_{i=1}^n u_i\eta_i\Big|\geq \frac{1}{3}\Big)&=
\Pr\Big(\Big|\sum_{i=1}^n u_i\eta_i\Big|^2\geq 
\frac{1}{3}\Ex \Big|\sum_{i=1}^n u_i\eta_i\Big|^2\Big)
\\
&\geq
\big(1-\frac{1}{3}\big)^2
\frac{(\Ex |\sum_{i=1}^n u_i\eta_i|^2)^2}{\Ex |\sum_{i=1}^n u_i\eta_i|^4}
\geq \frac{4}{27},
\end{align*}
thus
\[
\Pr(|\langle u,X\rangle|\geq t)\geq \frac{4}{27}\Pr(3(|g_1|+1)\geq t)
\geq \frac{4}{27}\Pr\Big(|g_1|\geq \frac{t}{3}\Big).
\]
Using the simple estimate 
$2t\exp(-(2t)^2/2)/\sqrt{2\pi}\leq \Pr(|g|\geq t)\leq \exp(-t^2/2)$,
we immediately get (\ref{example_1}) for $t\geq 3$.
For $0\leq t\leq 3$ we have  
\begin{align*}
\Pr(|\langle u,X\rangle|\leq t)&\leq
\Pr\Big(9\Big|\sum_{i=1}^{n}u_i\eta_i\Big|\leq t\Big)\leq 
\frac{\sqrt{2}t}{9}\leq t\frac{\Pr(|g_1|\leq 3)}{3} 
\\
&\leq \Pr(|g_1|\leq t)=
\Pr(|\langle u,Y\rangle|\leq t),
\end{align*}
where to get the second inequality we used Ball's upper bound on
cube sections \cite{B}. Hence (\ref{example_1}) holds also for $t\in [0,3]$.

Thus $Y\prec_{\omega} Y$. However
$\Ex\|Y\|=\Ex\max_{i\leq n}|g_i|\geq \sqrt{\log n}/L$ and
$\Ex\|X\|\leq 9\Ex(|g_1|+1)\leq 18$. 

\medskip

{\bf Acknowledgments.} The author would like to thank prof. S.\ Kwapie\'n
for suggesting the method used in the proof of Theorem \ref{str_tails}.

\noindent
Institute of Mathematics\\
Warsaw University\\
Banacha 2\\
02-097 Warszawa\\
{\tt rlatala@mimuw.edu.pl}

\end{document}